\documentclass[11pt]{article}
\usepackage{amsfonts}
\newcommand{\lvec}{\left[\begin{array}{c}}
\newcommand{\rvec}{\end{array}\right]}
\newcommand{\lmat}{\left[\begin{array}{cc}}
\newcommand{\rmat}{\end{array}\right]}
\newtheorem{theorem}{Theorem}

\newtheorem{lemma}{Lemma}

\begin{document}

\title{Convolution-type stochastic Volterra equations with additive fractional Brownian motion in Hilbert space}

\author{Peter Caithamer \\ Department of Mathematics and Actuarial Science \\ Indiana University Northwest \\ 3400 Broadway \\ Gary, IN 46408 \\ pcaitham@iun.edu \and Anna Karczewska \\ Department of Mathematics \\ University of Zielona G\'ora \\ ul.\ Szafrana 4a, 65-246 \\
Zielona G\'ora, Poland \\
A.Karczewska@im.uz.zgora.pl}

\maketitle

\begin{quote}{\bf Abstract} We consider convolution-type stochastic Volterra 
equations with additive Hilbert-valued fractional Brownian motion, $0<H<1$. 
We find the weak solution to this stochastic Volterra equation, and study its 
stochastic integral part, the stochastic convolution, which we show to be 
mean-zero Gaussian. We develop an It\^o isometry for stochastic integrals with 
respect to a Hilbert-valued fractional Brownian motion, and use it to compute 
the covariance of the stochastic convolution. This formula, which uses 
fractional integrals and derivatives, generalizes the well-known formula from 
the case $H=1/2$.    
AMS 2000 Subject Classification: 60H20,45D05, 60H05, 60G15
\newline Key words and phrases: fractional Brownian motion, fractional integrals and derivatives, stochastic convolution, stochastic integral equation, weak solutions.
\end{quote}

\section{Introduction}

Fractional Brownian motion is very useful for modeling systems which exhibit 
memory. These systems occur in many fields including finance and economics, physics, metallurgy, telecommunications, etc. The stochastic calculus with respect to a real-valued fractional Brownian has recently been developed, see \cite{hu}, \cite{nualart} and the references therein.

Meanwhile the stochastic calculus with respect to a Brownian motion or a martingale in Hilbert space have been frequently considered and are summarized in \cite{dapratozabczyk}, \cite{greckschtudor}. Recently, several authors have considered stochastic differential equations with respect to fractional Brownian motion in Hilbert space. The authors of \cite{dmpd1} define the stochastic integral w.r.t.\ fBm in Hilbert space for $H>1/2$ via a standard decomposition into a series of real-valued fBms and then study properties of mild solutions to a linear differential equation with an additive noise. A nonlinear differential equation with a multiplicative noise is studied for $H>1/2$ in \cite{maslowskinualart} via a similarly defined stochastic integral, where existence, uniqueness, and regularity results for the solution are obtained. A slightly more general equation but with an additive noise is considered for any value of $H$ in \cite{anhgrecksch}, which obtains existence and uniqueness results. The authors of \cite{maslowskischmalfuss} consider the stationarity of solutions of linear and nonlinear differential equations with an additive noise with respect to fBm in Hilbert space for $H>1/2$. The authors of \cite{ttv} consider for both $H>1/2$ and $H<1/2$ a linear differential equation with additive noise and a heat equation with multiplicative noise. They prove existence, uniqueness, and regularity results and establish a Feynman-Kac formula for the latter. A linear equation with multiplicative noise is considered for $H>1/2$ in \cite{dmpd2}, where reesults are obtained on the existence and stability of its solution. In \cite{djpd} a stochastic integration more general than that considered in \cite{dmpd1} is developed for $H>1/2$ and is used to prove an infinite dimensional It\^o formula.  

Linear stochastic evolution equations with memory are an emerging area of
research with interesting mathematical questions and various important
applications. These equations have been treated by many authors, see e.g.\
\cite{BoFa,BoTu,CDP96,CDP97,CDPP,RSS1,RSS2} or, very recently 
\cite{KaLi1,KaLi2,KaLi3}. In \cite{CDP96,CDP97,CDPP}, stochastic Volterra 
equations are studied in connection with viscoelasticity and heat conduction in materials with memory. The paper \cite{CDP96} due to Cl\'ement and Da Prato is particularly significant because the authors have extended the semigroup napproach, usually applied to stochastic differential equations, to a class of the Volterra equations. In the majority of the papers mentioned above, the resolvent approach is used. This approach allows one to answer the fundamental questions related to stochastic Volterra equations in Hilbert space like existence and regularity of mild, weak and even, in particular cases, strong solutions. Papers \cite{CDP96,CDP97,CDPP} provide regularity of stochastic convolutions corresponding to stochastic Volterra equations, particularly H\"olderianity of the corresponding trajectories. The paper \cite{CDP96} was the main source for the authors of \cite{BoFa,BoTu}, where mild solutions for the fractional and semilinear stochastic Volterra equations, respectively, have been studied. Papers \cite{BoTu} and \cite{RSS2} establish a large deviation principles for stochastic Volterra equations. The papers \cite{KaZa} and \cite{KaLi1} provide conditions  on the covariance of the spatially homogeneous noise term under which the solutions of the Volterra equations take values in a Sobolev space $H^\alpha$, particularly $L^2$. Volterra equations driven by cylindrical Wiener process have been considered in \cite{karczewska} and recently in \cite{KaLi2,KaLi3}. These papers generalize well-known result obtained for genuine Wiener noise or even classical results for differential equations. In \cite{KaLi2,KaLi3}, several versions of approximation theorems for resolvents are given, particularly for fractional Volterra equations. Those convergences, in some sense analogous to the Hille-Yosida theorem for semigroups, enable one to obtain strong solutions to a subclass of stochastic Volterra equations. 

In this paper we consider the convolution-type stochastic Volterra equation
\begin{equation}\label{sve}X(t) = X(0) + \int_0^t a(t-s)AX(s)ds + \int_0^t F(s) dB^H (s)\ ,\end{equation}
where $t\in \mathbb{R}_+$, $X(t)$ belongs to a Hilbert space $\mathcal{H}$, $a\in L^1_{\rm loc} (\mathbb{R}_+)$, $A$ is a closed linear operator mapping its domain $\mathcal{D}(A)$, which is assumed to be dense in $\mathcal{H}$, into $\mathcal{H}$,  $B^H$ is an $\mathcal{H}$-valued  fractional Brownian motion, $F(t)\in L(\mathcal{H}):= L(\mathcal{H},\mathcal{H})$ for any $t$. To the authors' knowledge this is the first paper considering stochastic Volterra equations in Hilbert space driven by a fractional Brownian motion.  

Section \ref{preliminaries} introduces background material on resolvents, fractional integrals and derivatives, stochastic integration w.r.t.\ Hilbert-valued fBm, and strong and weak solutions to stochastic Volterra equations, used in the rest of the paper. In Section \ref{stochconv} we find a weak solution the equation (\ref{sve}), find conditions under which it is unique, and show that its integral part, the stochastic convolution, is Gaussian. We develop an It\^o isometry formula for stochasic integrals w.r.t.\ Hilbert-valued fBm, and use it to compute the covariance of the stochastic convolution. Section \ref{examples} consider several examples to which the theory developed is applicable.

\section{Preliminaries}\label{preliminaries}

This section introduces background material on convolution-type stochastic Volterra equations, stochastic integration and integral equations in Hilbert space.

\subsection{Resolvents}\label{wpr} Equation (\ref{sve}) is a stochastic version of the deterministic Volterra equation  
\begin{equation}\label{dve}u(t) = u(0) + \int_0^t a(t-s)Au(s)ds\ .
\end{equation}
Equation (\ref{dve}) is said to be \emph{well-posed} if for each $u(0)\in \mathcal{D}(A)$ there exists a solution $u(t) \in \mathcal{D}(A)$ to (\ref{dve}) which is continuous in $t$ and $u_n(0)\to 0$ implies $u_n(t)\to 0$ uniformly on compact intervals. We will refer to equation (\ref{sve}) as \emph{well-posed} if equation (\ref{dve}) is.

The \emph{solution operator} or \emph{resolvent}, $S(t)$, is defined for $t\geq 0$ by 
\begin{equation}\label{resolvent}S(t)u(0)= u(t)\ .\end{equation}
Equation (\ref{dve}) is well-posed if and only if the solution operator is strongly continuous, commutes with $A$ and equation (\ref{resolvent}) solves equation (\ref{dve}). We note that $S(t)$ is generated by $A$ and $a$. More details, including conditions on $A$ and $a$ which guarantee the well-posedness of (\ref{dve}), may be found in \cite[Chapter 1]{pruss}, e.g.\ Theorem 1.1 p.\ 37. 

Formally substituting (\ref{resolvent}) into (\ref{dve}) we obtain the \emph{resolvent equation} 
\begin{equation}\label{resolventeq}S(t)u(0) = u(0) + \int_0^t a(t-s)AS(s)u(0)ds\ .\end{equation}
The resolvent is differentiable if, for example, $a\in BV_{\rm loc} (\mathbb{R}_+)$, \cite[p.\ 34]{pruss}. Differentiating (\ref{resolventeq}) yields
\begin{equation}\label{resolvderiv}\dot{S}(t)u(0) = \int_0^t \dot{a}(t-s)AS(s)u(0)ds + a(0)AS(t)u(0)\ .\end{equation}

\subsection{Fractional integrals and derivatives}

Define, as in \cite[p.\ 33]{samko}, for $f\in L^1(a,b)$ and $\alpha>0$ the left and right (respectively) Riemann-Liouville fractional integrals
\begin{equation}\label{lefti}(I_{a+}^\alpha f)(x)= \frac{1}{\Gamma (\alpha)}\int_a^x \frac{f(y)}{(x-y)^{1-\alpha}}dy\ ,\ \ \ x>a\ ,\end{equation}
and
\begin{equation}\label{righti}(I_{b-}^\alpha f)(x)= \frac{1}{\Gamma (\alpha)}\int_x^b \frac{f(y)}{(y-x)^{1-\alpha}}dy\ ,\ \ \ x<b\ . 
\end{equation}
The fractional derivatives are defined by
\begin{equation}\label{leftd}(D_{a+}^\alpha f)(x)=\frac{d^{\lceil\alpha\rceil}}{dx^{\lceil\alpha\rceil}} (I_{a+}^{\lceil\alpha\rceil-\alpha} f) (x)\ ,\end{equation}
and 
\begin{equation}\label{rightd}(D_{b-}^\alpha f)(x)=\frac{d^{\lceil\alpha\rceil}}{dx^{\lceil\alpha\rceil}} (I_{b-}^{\lceil\alpha\rceil-\alpha} f) (x)\ .\end{equation}

One may allow the order $\alpha$ of fractional integration/differentiation to be negative and by adopting the following notations $D^{-\alpha}=I^\alpha$ and $I^{-\alpha}=D^\alpha$. This allows one to write formulae compactly.

Let $\mathcal{K}$ be Hilbert space and $F \in L^1((a,b), \mathcal{K})$. Fractional integrals and derivatives of $F$ may be defined analogously with (\ref{lefti})-(\ref{rightd}) using a Bochner integral. Let $\{k_i\}_{i=1}^\infty$ be basis for $\mathcal{K}$, decompose $F(x) = \sum_{i=1}^\infty f_i (x) k_i$. Then for $D^\alpha$, the left/right integral/derivative, we have
\begin{equation}\label{fracddist}D^\alpha F(x) = \sum_{i=1}^\infty D^\alpha f_i (x) k_i\ .\end{equation} 
   
\subsection{Stochastic integration with respect to Hilbert-valued fBm}

A real-valued \emph{fractional Brownian motion (fBm)}, $b^H (t)$, with \emph{Hurst parameter},
$H\in(0,1)$, is a mean-zero Gaussian stochastic process with covariance function given by
\begin{equation}\label{fBmcovar}E[b^H (t)b^H (t')]=\frac{1}{2}(|t|^H + |t'|^H -|t-t'|^H)\ .\end{equation}
The increments of fBm are negatively correlated, independent, or positively correlated for
$H<1/2$, $H=1/2$, or $H>1/2$, respectively. The case $H=1/2$ corresponds to the usual
Brownian motion. Fractional Brownian motion is discussed in detail in \cite{mandelbrot}.
 
For $H>1/2$ define as in \cite{dhpd} 
\[\theta_H (s,t):= H(2H-1)|s-t|^{2H-2}= \frac{\partial^2}{\partial s\partial
t}\mbox{Cov}(b^H(s),b^H(t))\ .\]
For $T\in (0,\infty]$, denote by $L^2_{H}([0,T])$ the set of nonrandom functions for which 
\[\int_0^T\int_0^T| f(s)f(t)|\theta_H (s,t)dsdt<\infty\ .\]
It is proven in \cite{grippenberg} that for $f,g\in L^2_{H} ([0,T])$
\begin{equation}\label{fractionalito>}E\left[\int_0^T f(t)db^H_t\int_0^T g(t)db^H_t\right]= \int_0^T\int_0^T f(s)g(t)\theta_H (s,t)dsdt\ .\end{equation}

For $H<1/2$, $f,g \in L^2_H ([0,t]):= I_{t-}^{1/2-H} (L^2 ([0,t]))$ we have \cite{nualart}, \cite{perezabreu} 
\begin{eqnarray}\label{fractionalito<finite}E\left[\int_0^t f(s)db^H_s\int_0^t g(s)db^H_s\right]&=& \frac{2H \Gamma (3/2-H)\Gamma (H+1/2)}{\Gamma(2-2H)} \cdot \\
&&\int_0^t s^{1-2H}(D_{t-}^{1/2-H} s^{H-1/2}f) (s)(D_{t-}^{1/2-H} s^{H-1/2}g) (s)ds\ .\nonumber\end{eqnarray}

Note that this last formula also applies for $H>1/2$ and is consistent with (\ref{fractionalito>}). The fractional derivative is then interpreted as a fractional integral. 

Let $H\in (0,1)$, $B^H (t)$ be an $\mathcal{H}$-valued fractional Brownian motion, i.e.\ a mean zero Gaussian stochastic process with 
\begin{equation}\label{HfBm}E[B^H(t)\otimes B^H (t')]= \frac{1}{2}(|t|^H + |t'|^H -|t-t'|^H)\cdot \Lambda\ ,\end{equation}
where $\Lambda \in L^1 (\mathcal{H})$.

We define the (fractional) Brownian filtration $\mathcal{F}_t := \sigma (B^H (s); 0\leq s \leq t)$.

For each $t \in \mathbb{R}_+$ let $F(t) \in L(\mathcal{H})$, where for any $h \in \mathcal{H}$, $F(t)h \in L^2_H$. 

Noting the well-known decomposition
\begin{equation}\label{fBmdecomp}B^H (t) = \sum_{k=1}^{\infty} \sqrt{\lambda_k} h_k (x) b^H_k (t)\ ,\end{equation}
where  
$\{ b^H_k (t)\}_{k=1}^{\infty}$ is an independent and identically distributed (i.i.d.)\ sequence of real-valued fractional Brownian motions, and $\lambda_k$, $h_k$ are the eigenvalues and eigenfunctions of $\Lambda$, we may define 
\begin{equation}\label{Hintdef}\int_0^t F(s) dB^H (s) := \sum_{k=1}^{\infty} \sqrt{\lambda_k}\int_0^t F(s)h_k (x) db^H_k (t)\ .\end{equation}  

We note that since the integrands in (\ref{Hintdef}) are non-random, the stochastic integrals are both It\^o/Skorohod integrals and Stratonovich integrals, \cite{nualart}.

We also note that (\ref{Hintdef}) allows one to generate a stochastic Fubini theorem with respect to $\mathcal{H}$-valued fBm from a stochastic Fubini theorem for real-valued fBm.

\subsection{Strong vs.\ weak solutions}

Following \cite{dapratozabczyk} and \cite{sobczyk} we call $X(t)\in \mathcal{F}_t$ a \emph{strong solution} to equation (\ref{sve}) if $X(\cdot)$ has a version such that
\begin{itemize}
\item[(i)] for almost all $t\in [0,T]$,  $P(X(t)\in \mathcal{D}(A))=1$,
\item[(ii)] for any $t\in [0,T]$, $P\left( \int_0^t \|a(t-s)AX(s)\|_{\mathcal{H}} ds < \infty \right)=1$, and 
\item[(iii)] for any $t\in [0,T]$, equation (\ref{sve}) holds almost surely.
\end{itemize}

Similarly we call $X(t)\in \mathcal{F}_t$ a \emph{weak solution} to equation (\ref{sve}) if 
\begin{itemize}
\item[(i)] for any $t\in [0,T]$, $P\left( \int_0^t \|a(t-s)X(s)\|_{\mathcal{H}} ds < \infty \right)=1$, and 
\item[(ii)] for any $t\in [0,T]$ and any $\varphi \in \mathcal{D} (A^*)$ the following equation holds:
\begin{equation}\label{weak}\langle X(t), \varphi \rangle_\mathcal{H} = \langle X(0), \varphi \rangle_\mathcal{H} + \left\langle \int_0^t a(t-s) X(s) ds , A^*\varphi \right\rangle_\mathcal{H} + \left\langle \int_0^t F(s)B^H (s), \varphi \right\rangle_\mathcal{H}\ .\end{equation} 
\end{itemize}

We do not define a mild solution given by equation (\ref{weaksol}) below, since it is shown in the next section to be a weak solution.

Note that these defintions of strong and weak solutions follow those used in deterministic ODE/PDE theory and differ from those used for SDE theory by \cite{oksendal} and other authors. 

\section{Results on the stochastic convolution}\label{stochconv}

In these section we find the unique weak solution to (\ref{sve}), study its stochastic integral part (the stochastic convolution), and  gather properties of the stochastic convolution. We begin with a lemma.

\begin{lemma} If $a\in C^1$ and $X(t)$ is a weak solution to (\ref{sve}) and $\Phi \in C^1 ([0,T], \mathcal{D}(A^*))$, then
\begin{eqnarray} \langle X(t), \Phi (t) \rangle_\mathcal{H} &=& \langle X(0), \Phi (0) \rangle_\mathcal{H} + \int_0^t \langle \dot{a} \star X(s) + a(0) X (s), A^* \Phi (t) \rangle_\mathcal{H} ds \nonumber \\
&& + \int_0^t \langle \Phi (s), F(s)dB^H (s) \rangle_\mathcal{H} + \int_0^t\langle X(s), \dot{\Phi} (s) \rangle_\mathcal{H}ds\label{convoluted}\end{eqnarray}
\end{lemma}

{\bf Proof} Let $\{ \phi_i \}_{i=1}^\infty$ be an orthonormal basis for $\mathcal{H}$. Then $\Phi (t) = \sum_{i=1}^\infty f_i (t) \phi_i$. Focus on just one term of the sum and consider $\langle X(t), f(t)\phi \rangle_\mathcal{H} - \langle X(0), f(0)\phi \rangle_\mathcal{H}$
\begin{eqnarray}  &=& \int_0^t \frac{d}{du} \langle X(u), f(u)\phi \rangle_\mathcal{H} du \nonumber \\ 
&=& \int_0^t \frac{d}{du}\left\langle \int_0^u a(u-s) X(s) ds, \int_0^u \dot{f}(s) ds A^*\phi \right\rangle_\mathcal{H} du \nonumber\\
&&+\int_0^t \frac{d}{du}\left\langle \int_0^u F(s) dB^H (s), \int_0^u \dot{f}(s) ds \phi \right\rangle_\mathcal{H} du \nonumber \\ 
&=& \int_0^t \left\langle \int_0^u \dot{a}(u-s) X(s) ds + a(0)X(u), \int_0^u \dot{f}(s) ds A^*\phi \right\rangle_\mathcal{H} du \nonumber \\
&&+ \int_0^t \left\langle \int_0^u a(u-s) X(s) ds, \dot{f}(u) A^*\phi \right\rangle_\mathcal{H} du\nonumber \\
&&+\int_0^t \left\langle  F(u) dB^H (u), \int_0^u \dot{f}(s) ds \phi \right\rangle_\mathcal{H} du\ . 
\label{convlemma} \end{eqnarray}  
Summing (\ref{convlemma}) over $i$ completes the proof.

\begin{theorem}\label{theoweak} A weak solution of the convolution-type stochastic Volterra equation (\ref{sve}) is given by
\begin{equation}\label{weaksol}X (t) = S(t)X(0) + \int_0^t S (t-s) F(s) dB^H (s)\ .\end{equation}
If $a,S \in C^1$, then the solution is unique.
\end{theorem}

{\bf Proof.} First, we will show that any weak solution of (\ref{sve}) must satisfy (\ref{weaksol}). To do this, let $\phi \in \mathcal{D}(A^*)$, consider (\ref{convoluted}) with $\Phi (s) = S^* (t-s) \phi $ yielding
\begin{eqnarray}\label{unique}
\langle X(t), S^* (0) \phi \rangle_\mathcal{H} &=& \langle X(0), S^* (t) \phi \rangle_\mathcal{H} + \int_0^t \langle \dot{a} \star X(s) + a(0) X (s), A^* S^* (t-s) \phi \rangle_\mathcal{H} ds \nonumber \\
&& + \int_0^t \langle S^* (t-s) \phi, F(s)dB^H (s) \rangle_\mathcal{H} + \int_0^t\langle X(s), \frac{d}{ds}(S^* (t-s) \phi) \rangle_\mathcal{H}ds\ . 
\end{eqnarray}  
The sum of the first and third integral is zero by the weak form of (\ref{resolvderiv}), yielding the result.

To prove that (\ref{weaksol}) satisfies (\ref{weak}), we substitute (\ref{weaksol}) in (\ref{weak}) yielding 
\begin{eqnarray}&&\label{sf1}\left\langle S(t)X(0) + \int_0^t S (t-s) F(s) dB^H (s), \varphi \right\rangle_\mathcal{H} = \langle X(0), \varphi \rangle_\mathcal{H} \\
&&+ \left\langle \int_0^t a(t-s) \left[S(s)X(0) + \int_0^s S (s-u) F(u) dB^H (u)\right] ds , A^*\varphi \right\rangle_\mathcal{H} + \left\langle \int_0^t F(s)B^H (s), \varphi \right\rangle_\mathcal{H}\ .\nonumber
\end{eqnarray} 
We focus on second-to-last term in equation (\ref{sf1}), which we consider in two parts. First, 
\begin{equation} \left\langle \int_0^t a(t-s) S(s)X(0) ds , A^*\varphi \right\rangle_\mathcal{H} =  \left\langle S(t)X(0) , \varphi \right\rangle_\mathcal{H} - \left\langle X(0) , \varphi \right\rangle_\mathcal{H}\ .\end{equation}

Second, $\left\langle \int_0^t a(t-s) \int_0^s S (s-u) F(u) dB^H (u) ds , A^*\varphi \right\rangle_\mathcal{H}$
\begin{eqnarray*} &=& \left\langle \int_0^t \int_u^t a(t-s) S(s-u) F(u) ds dB^H (u) , A^*\varphi \right\rangle_\mathcal{H} \\
&=& \left\langle \int_0^t \int_0^{t-u} a(t-u-v) S(v) dv F(u) dB^H (u) , A^*\varphi \right\rangle_\mathcal{H}\\
&=& \left\langle \int_0^t (a \star S) (t-u) F(u) dB^H (u) , A^*\varphi \right\rangle_\mathcal{H}\\
&=& \left\langle \int_0^t [S (t-u) -I] F(u) dB^H (u) , \varphi \right\rangle_\mathcal{H}\\
&=& \left\langle \int_0^t S (t-u) F(u) dB^H (u) , \varphi \right\rangle_\mathcal{H} - \left\langle \int_0^t F(u) dB^H (u) , \varphi \right\rangle_\mathcal{H}\ ,
\end{eqnarray*}
where the first equality follows from a stochastic Fubini theorem, the second from a change of variables, the third from the definition of convolution, the fourth from the resolvent equation (\ref{resolventeq}). This yields the result.

We now focus on the stochastic integral part of (\ref{weaksol}) and define the stochastic convolution associated with the equation (\ref{sve}) as
\begin{equation}\label{stochconvol}X_{(A,a)}^F (t) = X_S^F (t) := \int_0^t S (t-s) F(s) dB^H (s)\ ,\end{equation}
where $S(t)$ is the resolvent obtain from $A$ and $a$. In order to study its properties we need the following version of the It\^o isometry for stochastic integrals with respect to $\mathcal{H}$-valued fBm with determnistic integrands.

\begin{theorem}\label{itoisom}
Let $B^H$ be an $\mathcal{H}$-valued fractional Brownian motion for $H\in (0,1)$. Let $F, G \in L^2_H ([0,t], L^2(\mathcal{H}))$. Then
\begin{eqnarray}\label{fractionalito}&&E\left[\int_0^t F(s)dB^H (s) \otimes\int_0^t G(s)dB^H (s)\right]= \\
&&\quad c(H) \int_0^t s^{1-2H}(D_{t-}^{1/2-H} s^{H-1/2}F) (s)\Lambda (D_{t-}^{1/2-H} s^{H-1/2}G)^* (s) ds\ ,\nonumber \end{eqnarray}
where $c(H) = (2H \Gamma (3/2-H)\Gamma (H+1/2))/ (\Gamma(2-2H))$.   
\end{theorem}

{\bf Proof.} Let $\{\lambda_i\}$ and $\{h_i\}$ be the eigenvalues and eigenvectors of $\Lambda$, the spatial covariance operator of $B^H (t)$. Noting that $L(\mathcal{H})$ may be identified with $\mathcal{H} \otimes \mathcal{H}$, decompose also $F(t) = \sum_{i,j=1}^{\infty} f_{i,j} (t) h_i \otimes h_j$. Using also (\ref{fBmdecomp}), we see 
\begin{eqnarray*}&&E\left[\int_0^t F(s)dB^H (s) \otimes \int_0^t G(s)dB^H (s)\right]\\ 
&=& E\left[\left(\sum_{k=1}^{\infty}\int_0^t \left( \sum_{i,j=1}^{\infty} f_{i,j} (s) h_i \otimes h_j\right)  \sqrt{\lambda_k} h_k db^H_k (s)\right)\right.\\
&& \otimes \left. \left(\sum_{k=1}^{\infty}\int_0^t \left(\sum_{i',j'=1}^{\infty} g_{i',j'} (s) h_{i'} \otimes h_{j'}\right) \sqrt{\lambda_k} h_{k} db^H_{k} (s)\right)\right] \\
&=& E\left[\left(\sum_{i,j=1}^{\infty} \sqrt{\lambda_j} \int_0^t f_{i,j} (s) h_i db^H_j (s)\right) \otimes \left(\sum_{i',j'=1}^{\infty} \sqrt{\lambda_{j'} } \int_0^t  g_{i',j'} (s) h_{i'} db^H_{j'} (s)\right)\right] \\
&=& c(H)\sum_{i,j,i'=1}^{\infty} \lambda_j h_{i} \otimes h_{i'} \left( \int_0^t s^{1-2H} (D^{1/2-H}_{t-} s^{H-1/2} f_{i,j}) (s) (D^{1/2-H}_{t-} s^{H-1/2} g_{i',j}) (s) ds\right)\ .
\end{eqnarray*}
On the other hand 
\begin{eqnarray*}&&\int_0^t s^{1-2H}(D_{t-}^{1/2-H} s^{H-1/2}F) (s)\Lambda (D_{t-}^{1/2-H} s^{H-1/2}G)^* (s)  ds \\
&=&\int_0^t s^{1-2H}\left(D_{t-}^{1/2-H} s^{H-1/2} \left(\sum_{i,j=1}^{\infty} f_{i,j} h_i \otimes h_j\right)\right) (s)\cdot \left(\sum_{k=1}^\infty \lambda_k h_k\otimes h_k\right) \\
&& \left(D_{t-}^{1/2-H} s^{H-1/2} \left(\sum_{i',j'=1}^{\infty} g_{i',j'}  h_{j'} \otimes h_{i'}\right)\right) (s) ds \\
&=& \int_0^t s^{1-2H}\left(\sum_{i,j=1}^{\infty} \lambda_j \left(D_{t-}^{1/2-H} s^{H-1/2} f_{i,j}\right) (s) h_i \otimes h_j\right) \left(\sum_{i',j'=1}^{\infty}\left(D_{t-}^{1/2-H} s^{H-1/2}  g_{i',j'}\right) (s) h_{j'} \otimes h_{i'}\right) ds \\
&=& \int_0^t s^{1-2H} \sum_{i,j,i'=1}^{\infty} (\lambda_j  (D_{t-}^{1/2-H} s^{H-1/2} f_{i,j}) (s) (D_{t-}^{1/2-H} s^{H-1/2} g_{i',j}) (s) ) h_{i} \otimes h_{i'} ds \ , 
\end{eqnarray*}
completing the proof.

\begin{theorem} If $S(t-s) F(s) \in L^2_H ([0,t], L^2(\mathcal{H}))$
then $ X_S^F (t) \in \mathcal{F}_t $ and 
\[X_S^F (t) \sim N\left(0, \int_0^t s^{1-2H} (D_{t-}^{1/2-H}s^{H-1/2}S(t-s)F(s))(s) \Lambda (D_{t-}^{1/2-H}s^{H-1/2}F^* (s) S^* (t-s))(s) ds \right)\ .\]
\end{theorem}

{\bf Proof.} Since the integrand the stochastic convolution $ X_S^F (t)$ is deterministic, it is immediate from the defintion of stochastic integration that $ X_S^F (t)$ is adapted and mean-zero Gaussian. To compute the covariance we simply apply Theorem \ref{itoisom}. We need only note that
\[(D_{t-}^{1/2-H}s^{H-1/2}F^* (s) S^* (t-s))(s)=(D_{t-}^{1/2-H}s^{H-1/2} S (t-s)F(s))^*(s)\ .\] 

We note that by using indicator functions one may compute $E[ X_S^F (t) \otimes X_S^F (t')]$. 

\section{Examples}\label{examples}

In this section we consider examples to which the results of this paper apply.

\subsection{Fractional integrodifferential equations}
Consider (\ref{dve}) where $\mathcal{H}$ is a space of functions of $x$, $Ah= (I_{c+}^\beta h)(x)$, where $\beta \in \mathbb{R}$, and $a(t-s)= (t-s)^{\alpha-1}/\Gamma (\alpha)$ for $\alpha >0$.  Then we have the following equation
\begin{equation}\label{dfrac}X(t)-X(0)= (I_{0+}^{\alpha} (I_{c+}^{\beta} X)(x)) (t) \ .\end{equation}
Since $\alpha>0$ then $a\in L^1_{\rm loc} (\mathbb{R}_+)$. But in order for the resolvent to be differentiable we need to require that $\alpha \geq 1$, so that $a\in BV_{\rm loc} (\mathbb{R}_+)$. Since the fractional dervative operator (\ref{leftd}) is the left inverse of the fractional integral operator (\ref{lefti}) we have (assuming that $X(0)=0$)
\[(D_{0+}^\alpha X)(t)= (I_{c+}^{\beta} X)(x)\ .\]
We thus consider the following stochastic equation
\begin{equation}\label{sfrac}X(t)-X(0)= (I_{0+}^{\alpha} (I_{c+}^{\beta} X)(x)) (t) + B^H (t)\ .
\end{equation}

We note that for $\alpha = 1$ (resp.)\ $2$ and $\beta= -2$ this is the stochastic heat (resp.)\ wave equation with an additive fractional noise.

We do not the specify the space $\mathcal{H}$, since we consider only weak solutions. 


\subsection{Parabolic Volterra equations with regular kernels} 
 
The equivalence of weak and mild solutions given in Theorem \ref{theoweak} is applicable for the whole class of the convolution type stochastic Volterra equations given by (\ref{sve}) (cf. \cite{karczewska}).

Now, let $a\in L^1_\mathrm{loc}(\mathbb{R}_+)$ be of subexponential growth, that is, $\int_0^{\infty}|a(t)|e^{-\epsilon t}dt< \infty$ for each $\epsilon >0$. The function $a$ is called $k$-\emph{regular}, where $k\in \mathbb{N}$, if there is a constant $c>0$ such that 
$ |\lambda^n \widehat{a}^{(n)}(\lambda)|\le c|\widehat{a}(\lambda)|$ for all 
$\mathrm{Re}\,\lambda >0$, $0\le n\le k$. Here $\widehat{a}$ denotes the Laplace transform of $a$.

Recall that (\ref{dve}) is called \emph{parabolic} if $\widehat{a}(\lambda)\ne 0$ for $\mathrm{Re}\,\lambda >0$, $\frac{1}{\widehat{a}(\lambda)}$ belongs to $\varrho(A)$, the resolvent set of the operator $A$, and
there is a constant $M>0$ such that $||(I-\widehat{a}(\lambda)A)^{-1}||\le $ for $\mathrm{Re}\,\lambda >0$. 

By \cite[Theorem 3.1]{pruss}, if (\ref{dve}) is parabolic and $a$ is 2-regular, then $S\in C^1((0,\infty);B(\mathcal{H}))$, where $B(\mathcal{H})$ denotes the space of all bounded linear operators on $\mathcal{H}$. Hence, if we restrict our considerations to the parabolic Volterra equations with the 2-regular kernel function $a$, then by Theorem \ref{theoweak} the weak solution is unique as well.

\end{document}